\documentclass{amsart}

\usepackage{amsmath}
\usepackage{amssymb}
\usepackage{amsthm}
\usepackage{mathtools}
\usepackage{graphicx}
\usepackage[nomessages]{fp}

\newtheorem{theorem}{Theorem}
\newtheorem{lemma}{Lemma}

\newcommand{\etal}{\emph{et al.}}
\newcommand{\abs}[1]{\lvert#1\rvert}
\newcommand{\pknk}{p_{k}(25n + 24 - k) \ \equiv \ 0 \pmod{5}}
\newcommand{\pkn}[1]{\FPeval\value{round(24 - #1,0)} p_{#1}(25n + \value) \ \equiv \ 0 \pmod{5}}
\newcommand{\cf}{\left(q; q\right)_{\infty}}

\newcommand{\cfq}[2]{\left(q^{#1}; q^{#2}\right)_{\infty}}
\newcommand{\cfl}[2]{\left(-q^{#1}; q^{#2}\right)_{\infty}}
\newcommand{\cfd}{\left(-q; q^{2}\right)_{\infty}}
\newcommand{\cfk}{\left(q; q^{2}\right)_{\infty}}

\begin{document}

\title{Some Observations on Modulo 5 Congruences for 2-Color Partitions}

\author{Suparno Ghoshal}
\address{Heritage Institute of Technology, Kolkata, India}
\email{ghoshalsuparno1331@gmail.com}

\author{Sourav Sen Gupta}
\address{Nanyang Technological University, Singapore}
\email{sg.sourav@gmail.com}

\date{}

\begin{abstract}
The 2-color partitions may be considered as an extension of regular partitions of a natural number $n$, with $p_{k}(n)$ defined as the number of 2-colored partitions of $n$ where one of the 2 colors appears only in parts that are multiples of $k$. In this paper, we record the complete characterization of the modulo 5 congruence relation $p_{k}(25n + 24 - k) \equiv 0 \pmod{5}$ for $k \in \{1, 2, \ldots, 24\}$, in connection with the 2-color partition function $p_k(n)$, providing references to existing results for $k \in \{1, 2, 3, 4, 7, 8, 17\}$, simple proofs for $k \in \{5, 10, 15, 20\}$ for the sake of completeness, and counter-examples in all the remaining cases. We also propose an alternative proof in the case of $k = 4$, without using the Rogers-Ramanujan ratio, thereby making the proof considerably simpler compared to the proof by Ahmed, Baruah and Ghosh Dastidar (JNT 2015).
\end{abstract}

\maketitle



\section{Introduction}
\label{sec:intro}

A \emph{partition} of a natural number $n$ is defined as a finite sequence of non-increasing positive integers $\{x_1, x_2, \ldots, x_s\}$ such that $n = x_1 + x_2 + \cdots + x_s$. The number of all such distinct partitions of $n$ is denoted by $p(n)$. It is quite well-known (and easy to prove) that the generating function of $p(n)$ is 
\[ \sum_{n=0}^\infty p(n) \, q^{n} 
\ = \ \prod_{\ell = 1}^{\infty} \left(\frac{1}{1 - q^{\ell}}\right) 
\ = \ \prod_{\ell = 0}^{\infty} \frac{1}{(1 - q \cdot q^{\ell})} \ = \ \frac{1}{\cf}, \]
where, for complex numbers $a$ and $\abs{q} < 1$, the $q$-shifted factorial is defined by
\[ (a;q)_\infty \ = \ \prod_{\ell = 0}^\infty (1 - aq^{\ell}). \]

The following identity for $p(n)$, discovered by Ramanujan~\cite{Ramanujan1927}, has been considered by stalwarts like G.~H.~Hardy and P.~A.~MacMahon as the ``most beautiful identity'' for the partition function (refer to~\cite{Hirschhorn2011} for a simple proof):
\begin{equation}
\label{eqn:beauty}
\sum_{n=0}^\infty p(5n+4) \, q^{n} \ = \ 5 \ \frac{\cfq{5}{5}^{5}}{\cf^{6}}
\end{equation}
This directly implies one of the most popular congruence relations involving the partition function $p(n)$, discovered by Ramanujan~\cite{Ramanujan1919} in 1919:
\[ p(5n+4) \ \equiv \ 0 \pmod{5}. \]

\subsection{2-color partitions}
\label{sec:twocolor}

2-color partitions may be considered as an extension of regular partitions of a natural number $n$, with $p_{k}(n)$ defined as the number of 2-color partitions of $n$ where one of the colors appears only in parts that are multiples of $k$. The generating function of $p_{k}(n)$, for natural numbers $n$ and $k$, as per~\cite{AhmedBaruahGD2015}, is
\begin{equation}
\label{eqn:genfun}
\sum_{n=0}^\infty p_{k}(n) \, q^{n} 
\ = \ \prod_{\ell = 1}^{\infty} \left(\frac{1}{1 - q^{\ell}}\right) \left(\frac{1}{1 - q^{k\ell}}\right) 
\ = \ \frac{1}{\cf \cfq{k}{k}}.
\end{equation}
Analogous to the congruence results of $p(n)$, recent papers~\cite{ChenLin2009,Chan2010,AhmedBaruahGD2015,Chern2016} have introduced several congruence results for the 2-color partition function $p_{k}(n)$. One of the most general forms for congruence in this direction is the following:
\begin{equation}
\label{eqn:maincong}
\pknk \qquad \text{for certain } k \in \{1, 2, \ldots, 24\}.
\end{equation}

\medskip

{\flushleft\textbf{Complete characterization} of this congruence for \emph{specific values} of $k$ is as follows:}
\begin{itemize}
\setlength{\itemsep}{5pt}

\item The cases $k = 6, 9, 11, 12, 13, 14, 16, 18, 19, 21, 22, 23, 24$ are \emph{false}, in fact, they all fail the congruence for $n = 0$, as $p_{6}(18) = 487$, $p_{9}(15) = 187$, $p_{11}(13) = 103$, $p_{12}(12) = 78$, $p_{13}(11) = 56$, $p_{14}(10) = 42$, $p_{16}(8) = 22$, $p_{18}(6) = 11$, $p_{19}(5) = 7$, $p_{21}(3) = 3$, $p_{22}(2) = 2$, $p_{23}(1) = 1$, and $p_{24}(0) = 1$.

\item The case $k = 1$ was proved by Baruah \etal~\cite{BaruahSarmah2013} in 2013. It is trivial.

\item The cases $k = 5, 10, 15, 20$ are easy to prove, as shown in Theorem~\ref{thm:k5cong}.

\item The cases $k = 3, 4$ were proved by Ahmed \etal~\cite{AhmedBaruahGD2015} in 2015. In this paper, we provide an alternative proof for $k = 4$ that is simpler than that of~\cite{AhmedBaruahGD2015}.

\item The cases $k = 7, 8, 17$ were conjectured by Ahmed \etal~\cite{AhmedBaruahGD2015} in 2015, and subsequently proved by Chern~\cite{Chern2016} in 2016.

\item The case $k = 2$ appears to be most frequent in the literature. It was first proved by Chen and Lin~\cite{ChenLin2009} in 2009, using modular forms. In 2010, Chan and Toh~\cite{ChanToh2010}, and in 2011, Xiong~\cite{xiong2011}, independently proved that for $\alpha \geq 2$,
\[ p_2(5^{\alpha}n + \delta_{\alpha}) \ \equiv \ 0 \pmod{5^{\lfloor \alpha / 2 \rfloor}}, \]
where $\delta_{\alpha}$ is the multiplicative inverse of $8$ modulo $5^{\alpha}$. Specifically for $\alpha = 2$, this reduces to $\pkn{2}$, that is, the case $k = 2$ for our 2-color partitions. In 2015, Ahmed \etal~\cite{AhmedBaruahGD2015} proposed an alternative proof to the same using known identities for the Ramanujan's theta function and the Rogers-Ramanujan continued fraction. In 2017, Hirschhorn described an (yet unpublished) algorithm to find formulas for $\sum_{n=0}^{\infty} p_2(5^{\alpha}n + \delta_{\alpha}) \, q^n$, where the case $\alpha = 2$ naturally includes $\pkn{2}$.

\end{itemize}

\subsection{Motivation and Contribution}
\label{sec:motivation}

In this paper, we record the complete characterization of the congruence relation $\pknk$, for all values of $k \in \{1, 2, \ldots, 24\}$, of the 2-color partition function $p_k(n)$, originally proposed by Ahmed \etal~\cite{AhmedBaruahGD2015} in 2015. In this direction, we have referred to the existing results (as in the current literature) in the previous section.

The main technical result of this paper is to provide an alternative proof of the congruence relation for the case $k = 4$, which is much simpler compared to the same by Ahmed \etal~\cite{AhmedBaruahGD2015}. Note that Ahmed \etal~\cite{AhmedBaruahGD2015} exploited known identities related to the Rogers-Ramanujan ratio in the course of their proof, while we provide a straight-forward approach avoiding the use of this ratio. 

The main result we prove in this paper, for the case of $k = 4$, is as follows:
\begin{theorem}
\label{thm:k4cong}
For any non-negative integer $n$, 
\begin{equation}
\label{eqn:k4cong}
\pkn{4}.
\end{equation}
\end{theorem}

We also prove the cases of $k \in \{5, 10, 15, 20\}$, for the sake of completeness:
\begin{theorem}
\label{thm:k5cong}
If $k \in \{5, 10, 15, 20\}$, then for any non-negative integer $n$, 
\begin{equation}
\label{eqn:k5cong}
\pknk.
\end{equation}
\end{theorem}



\section{Preliminaries}
\label{sec:prelim}

In this paper, we use Ramanujan's general theta function and some of it's special cases as our main analytic tools to deconstruct the partition functions. Ramanujan's general theta function $f(a,b)$ is defined as follows (see page 34 of~\cite{NotebooksIII} for reference):
\[ f(a,b) \coloneqq \sum_{k=-\infty}^{\infty} a^{\frac{k(k+1)}{2}} \, b^{\frac{k(k-1)}{2}}, \qquad a, b \in \mathbb{C}, \ \lvert ab \rvert < 1. \]
The most important special case of $f(a,b)$ that we require for this paper is:
\begin{align}
\label{eqn:phi}
\phi(q) &= f(q,q) = \sum_{k=-\infty}^{\infty} q^{k^2} = \cfd^{2} \cfq{2}{2} = \frac{\cfq{2}{2}^{5}}{\cf^{2}\cfq{4}{4}^{2}} \end{align}

In addition, we also require the following identities involving the partition function, and the functions $\phi(\cdot)$ and $f(\cdot,\cdot)$, as mentioned above. The reader may refer to Ramanujan Notebooks -- Part III~\cite{NotebooksIII} for details of these results.

\begin{lemma}[Jacobi's Identity --- see page 39, Entry 24 (ii) of~\cite{NotebooksIII} for reference]
\label{lem:jacobi}
\[ \cf^{3} \ = \ \sum_{n=0}^\infty(-1)^{n} (2n+1) \, q^{\frac{n(n+1)}{2}} \]
\end{lemma}

\begin{lemma}[5-dissection of $\phi$ --- see page 49, Corollary (i) of~\cite{NotebooksIII} for reference]
\label{lem:phifive}
\[ \phi(q) \ = \ \phi(q^{25}) + 2 \, q \, f(q^{15}, q^{35}) + 2 \, q^{4} \, f(q^{5}, q^{45}) \]
\end{lemma}

\begin{lemma}[5-degree $\phi-f$ identity --- see page 262, Entry 10 (iv) of~\cite{NotebooksIII} for reference]
\label{lem:phidentity}
\[ \phi^{2}(q)-\phi^{2}(q^{5}) = 4 \, q \, f(q^{3}, q^{7}) \, f(q, q^{9}) \]
\end{lemma}



\section{Proof of Theorem~\ref{thm:k4cong}}
\label{sec:k4proof}

To prove $\pkn{4}$ for any non-negative integer $n$, we study the generating function of $p_4(n)$. From the identity for the generating function of $p_k(n)$, as in equation~\eqref{eqn:genfun}, we obtain the generating function of $p_4(n)$ as follows:
\[ \sum_{n=0}^\infty p_{4}(n) \, q^{n} \ = \ \frac{1}{\cf \cfq{4}{4}} \ = \ \frac{\cf^{5} \cfq{4}{4}^{5}}{\cf^{6} \cfq{4}{4}^{6}} \]
Multiply and divide the equation above by $\cfq{2}{2}^{15}$, and use~\eqref{eqn:phi} to obtain
\[ \sum_{n=0}^\infty p_{4}(n) \, q^{n} = \frac{\cf^{5} \cfq{4}{4}^{5}}{\cfq{2}{2}^{15}} \ \frac{\cfq{2}{2}^{15}}{\cf^6 \cfq{4}{4}^6} = \frac{\cf^{5} \cfq{4}{4}^{5}}{\cfq{2}{2}^{15}}  \, \phi^{3}(q). \]
Simple application of binomial expansion, with the equation above, produces the following equivalence relation (all equivalence relations in this proof are modulo 5):
\begin{align*}
& \sum_{n=0}^\infty p_{4}(n) \, q^{n} \ \equiv \ \frac{\cfq{5}{5}\cfq{20}{20}}{\cfq{10}{10}^{3}}  \, \phi^{3}(q) \\
& \equiv \ \frac{\cfq{5}{5}\cfq{20}{20}}{\cfq{10}{10}^{3}} \, (\phi(q^{25})+2qf(q^{15},q^{35})+2q^{4}f(q^{5},q^{45}))^{3}
\end{align*}
where we have also used the 5-dissection of $\phi(q)$, as in Lemma~\ref{lem:phifive}.

Extract \emph{only} the terms in which the power of $q$ is a multiple of $5$ to get
\[ \sum_{n=0}^\infty p_{4}(5n) \, q^{5n} \equiv \frac{\cfq{5}{5}\cfq{20}{20}}{\cfq{10}{10}^{3}} \, \phi(q^{25}) \left(\phi^{2}(q^{25})+4 q^{5} f(q^{15}, q^{35}) f(q^{5}, q^{45}) \right) \]
and thereafter replace $q^{5}$ by $q$ in the extracted terms to obtain
\[ \sum_{n=0}^\infty p_{4}(5n) \, q^{n} \ \equiv \ \frac{\cf\cfq{4}{4}}{\cfq{2}{2}^{3}} \, \phi(q^{5}) \left( \phi^{2}(q^{5})+4qf(q^{3},q^{7})f(q,q^{9}) \right). \]

Using the 5-degree $\phi-f$ identity (Lemma~\ref{lem:phidentity}), a series of manipulations involving the use of~\eqref{eqn:phi}, twice, and some simple applications of binomial expansion, we get
\begin{align*}
\sum_{n=0}^\infty p_{4}(5n) \, q^{n} \ & \equiv \ \frac{\cf\cfq{4}{4}}{\cfq{2}{2}^{3}} \, \phi(q^{5}) \, \phi^{2}(q)\\
& \equiv \ \phi(q^{5}) \, \frac{\cf\cfq{4}{4}}{\cfq{2}{2}^{3}} \, \frac{\cfq{2}{2}^{10}}{\cf^{4}\cfq{4}{4}^{4}} \\
& \equiv \ \phi(q^{5}) \, \frac{\cfq{2}{2}^{7}}{\cf^{3}\cfq{4}{4}^{3}} \\
& \equiv \ \phi(q^{5}) \, \frac{\cfq{2}{2}^{5}}{\cf^{5}\cfq{4}{4}^{5}} \, \cf^{2}\cfq{2}{2}^{2}\cfq{4}{4}^{2} \\
& \equiv \ \phi(q^{5}) \, \frac{\cfq{10}{10}}{\cfq{5}{5}\cfq{20}{20}} \, \cf^{2}\cfq{2}{2}^{2}\cfq{4}{4}^{2} \\
& \equiv \ \frac{\cfq{10}{10}^{6}}{\cfq{5}{5}^{3}\cfq{20}{20}^{3}} \, \cf^{2}\cfq{2}{2}^{2}\cfq{4}{4}^{2}.
\end{align*}

We may write $\cf = \cfk\cfq{2}{2}$, and similarly write $\cfq{5}{5} = \cfq{5}{10}\cfq{10}{10}$ to get
\[ \sum_{n=0}^\infty p_{4}(5n) \, q^{n} \ \equiv \ \frac{\cfq{10}{10}^{3}}{\cfq{5}{10}^{3}\cfq{20}{20}^{3}} \, \cfk^{2}\cfq{2}{2}^{4}\cfq{4}{4}^{2}. \]
Now replace $q$ by $-q$ in the above congruence to obtain
\begin{align*}
\sum_{n=0}^\infty (-1)^{n} \, p_{4}(5n) \, q^{n} \ & \equiv \ \frac{\cfq{10}{10}^{3}}{\cfl{5}{10}^{3}\cfq{20}{20}^{3}}  \, \cfd^{2}\cfq{2}{2}^{4}\cfq{4}{4}^{2}.
\end{align*}
Multiply the numerator and denominator by $\cfq{5}{10}^3$ and $\cfk^2$ to get
\[ \sum_{n=0}^\infty (-1)^{n} \, p_{4}(5n) \, q^{n} \ \equiv \ \frac{\cfq{5}{10}^{3}\cfq{10}{10}^{3}}{\cfq{10}{20}^{3}\cfq{20}{20}^{3}} \, \frac{\cfq{2}{4}^{2}\cfq{2}{2}^{4}\cfq{4}{4}^{2}}{\cfk^2}. \]
By simple counting of powers of $q$ in the binomial expansion of the terms, and using the even-odd power identity $\cf = \cfk \cfq{2}{2}$, we get
\begin{align*}
\sum_{n=0}^\infty (-1)^{n} \, p_{4}(5n) \, q^{n} \ & \equiv \ \frac{\cfq{5}{5}^{3}}{\cfq{10}{10}^{3}} \, \frac{\cfq{2}{2}^{8}} {\cf^{2}}\\
& \equiv \ \frac{\cfq{5}{5}^{3}}{\cfq{10}{10}^{3}} \, \frac{\cfq{2}{2}^{5}}{\cf^{5}} \, \cf^{3} \cfq{2}{2}^{3}.
\end{align*}
We know that $\cf^5 \equiv \cfq{5}{5} \pmod{5}$ and $\cfq{2}{2}^5 \equiv \cfq{10}{10} \pmod 5$, by simple binomial expansion of the terms. This reduces the above to
\begin{equation}
\label{eqn:simple}
B(q) \ := \ \sum_{n=0}^\infty (-1)^{n} \, p_{4}(5n) \, q^{n} \ \equiv \ \frac{\cfq{5}{5}^{2}}{\cfq{10}{10}^{2}} \, \cf^{3} \cfq{2}{2}^{3}.
\end{equation}

Expand $\cf^3$ and $\cfq{2}{2}^3$ using the Jacobi's Identity from Lemma~\ref{lem:jacobi} to obtain
\[ B(q) = \frac{\cfq{5}{5}^{2}}{\cfq{10}{10}^{2}} \left( \sum_{r=0}^\infty(-1)^{r} (2r+1) \, q^{\frac{r(r+1)}{2}} \right) \left( \sum_{s=0}^\infty(-1)^{s} (2s+1) \, q^{s(s+1)} \right). \]
Multiplying the summation terms in the above expression, we get
\begin{equation}
\label{eqn:simplecont}
B(q) = \frac{\cfq{5}{5}^{2}}{\cfq{10}{10}^{2}} \left( \sum_{r, s = 0}^\infty (-1)^{r+s} (2r+1) (2s+1) \, q^{\frac{r(r+1)}{2} + s(s+1)} \right).
\end{equation}

In the expression above, let us only consider the terms with power of $q$ of the form $5n + 4$. In all such terms, $\frac{r(r+1)}{2} + s(s+1) \equiv 4 \pmod{5}$. We know that $\frac{r(r+1)}{2}$ only assumes values $\{0, 1, 3\}$ modulo 5, and $s(s+1)$ only assumes values $\{0, 1, 2\}$ modulo 5. Thus, the only possible combination of $(r, s)$ that may result in the desired value $\frac{r(r+1)}{2} + s(s+1) \equiv 4 \pmod{5}$ is $(r = 2 \pmod{5}, s = 2 \pmod{5})$, and in this case, the coefficient turns out to be
\[ (-1)^{r+s} (2r+1) (2s+1) \ \equiv \ 0 \pmod{5} \]

Now, if we extract only the terms of the form $q^{5n + 4}$ from Equation~\eqref{eqn:simplecont}, replacing $B(q)$ by its series form $\sum_{n=0}^\infty (-1)^{n} \, p_{4}(5n) \, q^{n}$ on the left hand side, and use the result above on the coefficients of all such terms, we obtain
\[ \sum_{n=0}^\infty (-1)^{n} \, p_{4}(5(5n+4)) \, q^{5n + 4} = \sum_{n=0}^\infty (-1)^{n} \, p_{4}(25n+20) \, q^{5n + 4} \ \equiv \ 0 \pmod{5}. \]
It directly follows from the above that $p_{4}(25n+20) \equiv 0 \pmod{5}$ for all $n$. \qed



\section{Proof of Theorem~\ref{thm:k5cong}}
\label{sec:k5proof}

Note that to prove the desired result --- $p_k(25n + 24 - k) \equiv 0 \pmod{5}$ for $k \in \{5, 10, 15, 20\}$ and for any non-negative integer $n$, it is sufficient to prove that $p_{5\ell}(25n - 5\ell + 20 + 4) \equiv 0 \pmod{5}$ for $\ell \in \{1, 2, 3, 4\}$ and for any non-negative integer $n$. In turn, it is sufficient (in fact, it is stronger) to prove that $p_{5\ell}(5m + 4) \equiv 0 \pmod{5}$ for $\ell \in \{1, 2, 3, 4\}$ and for any non-negative integer $m$, where $5m = 5(5n - \ell + 4) \geq 0$ for $\ell \in \{1, 2, 3, 4\}$ and $n \geq 0$.

To prove $p_{5\ell}(5m + 4) \equiv 0 \pmod{5}$ for $\ell \in \{1, 2, 3, 4\}$ and for any non-negative integer $n$, we study the generating function of $p_{5\ell}(n)$. As per equation~\eqref{eqn:genfun}, the generating function of $p_{5\ell}(n)$ is as follows:
\[ \sum_{n=0}^\infty p_{5\ell}(n) \, q^{n} \ = \ \frac{1}{\cf \cfq{5\ell}{5\ell}} \]

We have to extract the terms of the form $q^{5m+4}$ from the above expression, where $\ell \in \{1, 2, 3, 4\}$. This is equivalent to extracting the terms of the form $q^{5m+4}$ from just the $1/\cf$ part in the above expression, as the $1/\cfq{5\ell}{5\ell}$ part contributes only to terms of the form $q^{5m}$ for any $\ell \geq 1$. 

Since we know that $1/\cf = \sum_{n=0}^\infty p(n) \, q^{n}$ and $p(5n+4) \equiv 0 \pmod{5}$ by Ramanujan~\cite{Ramanujan1927}, we can conclude that the terms of the form $q^{5m+4}$ extracted from the expression for $\sum_{n=0}^\infty p_{5\ell}(n) \, q^{n}$ will all have coefficients equal to zero modulo 5. Therefore, $p_{5\ell}(5m + 4) \equiv 0 \pmod{5}$ for $\ell \in \{1, 2, 3, 4\}$ and for any non-negative integer $m$. Hence the desired result. \qed



\section{Conclusion}
\label{sec:conclusion}
We record the complete characterization of the modulo 5 congruence relation $\pknk$ of the 2-color partition function $p_k(n)$, and propose an alternative proof in the case of $k = 4$, which is simpler than the proof by Ahmed, Baruah and Ghosh Dastidar~\cite{AhmedBaruahGD2015}. Note that the proof could be made even simpler if the individual terms are considered modulo 5 before multiplying them out to obtain Equation~\eqref{eqn:simplecont} in our proof of Theorem~\ref{thm:k4cong}. This approach is quite commonly used by Hirschhorn~\cite{Hirschhorn2017}, while multiplying the terms first, as we did while proving Theorem~\ref{thm:k4cong}, is a style generally attributed to Ramanujan. It would have also been nice to have alternative proofs for the cases $k = 7, 8, 17$, which are analytic and simpler compared to the modular form based proofs by Chern~\cite{Chern2016}.

\bigskip

{\flushleft\emph{Acknowledgement ---}} The authors would like to express their sincere gratitude to Professor Michael D. Hirschhorn for his constructive criticism and kind guidance throughout the course of the research project. His detailed comments on the initial drafts helped us a lot in improving the technical and editorial quality of this paper.



\bibliographystyle{amsplain} 
\bibliography{partitionref}

\end{document}